\magnification\magstep1

\font\sevenrm = cmr7

\centerline{\bf Complex Monge-Amp\`ere of a Maximum}


\bigskip

\centerline{E. Bedford and S. Ma`u\footnote{}{\sevenrm E.B. was supported in part by the NSF, and S.M. was supported
by a New Zealand Science \& Technology Post-Doctoral fellowship.} }

\bigskip

\bigskip

Pluri-subharmonic (psh) functions play a primary role in
pluri-potential theory.  They are closely related to the operator
$dd^c=2i\partial\bar\partial$ (with notation
$d=\partial+\bar\partial$ and $d^c=i(\bar\partial-\partial)$), which
serves as a generalization of the Laplacian from ${\bf C}$ to ${\bf
C}^{\rm dim}$ for ${\rm dim}>1$.    If $u$ is smooth of class $C^2$,
then for $1\le n\le {\rm dim}$, the coefficients of the exterior
power $(dd^cu)^n$ are given by the $n\times n$ sub-determinants of
the matrix $(\partial^2 u/\partial z_i\partial\bar z_j)$.  The top
exterior power corresponds to $n={\rm dim}$, and in this case we
have the determinant of the full matrix, which gives the complex
Monge-Amp\`ere operator. The extension of the (nonlinear) operator
$(dd^c)^n$ to non-smooth functions has been studied by several
authors (see, for instance, [B2]).  Here it will suffice to define
$(dd^c)^n$ on psh functions which are continuous.

For an open set  $E\subset{\bf C}^{\rm dim}$, the set of psh
functions on $E$ forms a cone which is closed under the operations
of addition and of taking finite maxima.   The relation between
$(dd^c)^n$ and the additive structure is given by the formula
$$(dd^c(u_1+u_2))^n=\sum_{n_1+n_2=n}{n!\over n_1!n_2!}\,(dd^cu_1)^{n_1}\wedge(dd^cu_2)^{n_2},$$
which holds in both the classical and generalized sense.  Here we
make the connection between $(dd^c)^n$ and the operation of  taking
a finite maximum.  A few examples of this are known already: if
$u=\max(0,\log(|z_1|^2+\cdots+|z_k|^2))$, then $(dd^cu)^k$ is a
multiple of surface measure on $\{|z_1|^2+\cdots+|z_k|^2=1\}$.  And
if $u=\max(0,x_1,\dots,x_k)$, then  $(dd^cu)^k$ is a multiple of
surface measure on $\{x_1=\cdots=x_k=0\}$. The case 
$(dd^c\max(u_1,u_2,u_3))^2$, $u_j$'s pluri-harmonic, is given in
[M].
 A different sort of formula for $(dd^c\max(u_1,u_2))^k$ is
given in [B1].  

More generally, let the functions  $u_1,\dots,u_m$ be smooth, and
let
$$u:=\max(u_1,\dots,u_m)$$
be their maximum.  It follows that $dd^cu_j$ and $dd^cu$ are locally
bounded below, and for the purpose of defining $(dd^c)^n$, $u_j$ and
$u$ may be treated as psh functions.   Taking the maximum stratifies
$E$ as follows: for each $J\subset\{1,\dots,m\}$  there is
$$E_J=E_J(u_1,\ldots,u_m):=\{z\in E:u(z)=u_j(z)\ \forall j\in J,\ u(z)>u_i(z)\ \forall i\notin J\}.$$
The sets $E_J$ form a partition of $E$, and as $J$ increases, the
subsequent $E_J$'s lie inside the ``boundaries'' of the previous
ones.  Since $(dd^cu)^n$ is representable by integration, we may
decompose it into a sum over the elements of the partition:
$$(dd^cu)^n=\sum_J (dd^cu)^n|_{E_J}.\eqno(1)$$
The first terms in (1) are easy to identify:  $E_j$ is open, and
thus $(dd^cu)^n|_{E_j}=(dd^cu_j)^n|_{E_j}$.  For the rest of the
terms, we use the following notation:   we write
$J=\{j_1,\dots,j_\ell\}$, $j_1<\cdots<j_\ell$, so $|J|=\ell$ is the
number of elements, and we define the forms
$$\delta_J^c=\delta_J^c(u_1,\dots,u_m):=d^c(u_{j_1}-u_{j_2})\wedge\cdots\wedge d^c(u_{j_{\ell-1}}-u_{j_\ell})\eqno(2)$$
$$\sigma_{J}^n=\sigma_{J}^n(dd^cu_1,\dots,dd^cu_m):=\sum_{\beta_1+\cdots+\beta_\ell=n}(dd^cu_{j_1})^{\beta_1}\wedge\cdots\wedge(dd^cu_{j_\ell})^{\beta_\ell}.\eqno(3)$$
If $E_J$ is smooth, we let $[E_J]$ denote the current of integration
over $E_J$, where we orient $E_J$ so that the current
$\delta^c_J\wedge [E_J]$ is positive (see Lemma 1).  In this paper
we identify the terms of the summation (1) as integrations on the
strata $E_J$: \proclaim Theorem 1.  Let $u_j\in C^3$, $1\le j\le m$,
be given, and set $u=\max(u_1,\dots,u_m)$.  If all of the sets $E_J$
are smooth, then
$$(dd^cu)^n=\sum_J \sigma_{J}^{n-|J|+1}\wedge \delta^c_J\wedge [E_J],\eqno(\star)$$
where the sum is taken over all $J\subset\{1,\dots,m\}$ with $1\le |J|\le n+1$.

A function that arises frequently is
$v=\log(|f_1|^2+\cdots+|f_N|^2)$, where the $f_i$'s are holomorphic.
Set $f=(f_1,\dots,f_N):{\bf C}^{\rm dim}\to{\bf C}^N$ and let
$\pi:{\bf C}^N-0\to{\bf P}^{N-1}$ be the projection.  The powers
$(dd^cv)^n$ may be determined on $\{v>-\infty\}$ using the fact that
$dd^c\log(|z_1|^2+\cdots+|z_N|^2)=\pi^*\omega_{FS}$, where
$\omega_{FS}$ denotes the Fubini-Study K\"ahler form on ${\bf
P}^{N-1}$.  That is, $(dd^cv)^n$ is the pullback of $\omega_{FS}^n$
under $(\pi\circ f)^*$.  In particular, $(dd^cv)^n=0$ if $n\ge N$.

Such functions and their maxima arise naturally with generalized
polyhedra.  For $1\le\alpha\le A$,  let
$p_{\alpha,1},\dots,p_{\alpha,N_\alpha}$ be polynomials, and let
${\rm deg}_\alpha=\max_{1\le i\le N_\alpha}({\rm deg}(p_{\alpha,
i}))$ be the maximum of the degrees.  Define
$$u_\alpha={1\over 2\,{\rm deg}_\alpha}\log\sum_{i=1}^{N_\alpha} |p_{\alpha,i}|^2  {\rm \ \ and\ }\  K=\{z\in{\bf C}^{\rm dim}:\max_{1\le \alpha\le A} (|p_{\alpha,1}|^2+\cdots+|p_{\alpha,N_{\alpha}}|^2)\le 1  \}.$$
A useful fact is: {\sl
$(dd^c(\max(u_{\alpha_1},\dots,u_{\alpha_\ell}))^N=0$ on the set
where the maximum is finite, whenever  $N\ge
N_{\alpha_1}+\cdots+N_{\alpha_\ell}$. }  This may be seen because
$(\star)$ involves sums of terms of the form
$(dd^cu_{a_1})^{\beta_1}\wedge\cdots\wedge(dd^cu_{a_k})^{\beta_k}$
with $\beta_1+\cdots+\beta_k=N+1-k$, which means that for some $i$,
$\beta_i\ge N_{\alpha_i}$ and thus $(dd^cu_{\alpha_i})^{\beta_i}=0$.

\proclaim Theorem 2.   Let us set $u:=\max(0,u_1,\dots,u_A)$.
Suppose that $u\ge\log^+|z|-C$ and that for every $J$ such that
$E_J-K\ne\emptyset$ we have $\sum_{j\in J}N_j\le {\rm dim} $.  Then
$u$ is the psh Green function for $K$; and if the sets $E_J$ are
smooth at points of $K$, then the equilibrium measure
$\mu_K:=({1\over 2\pi}dd^cu)^{\rm dim}$ is given by the formula
($\star$).

\noindent{\it Proof. }  We will show first that $u$ is the psh Green
function of $K$.  It is evident that $u$ is continuous and bounded
above by $\log^+|z|+C$.  Since we also have the lower bound
$\log^+|z|-C$, it will suffice to show that $(dd^cu)^{\rm dim}$
vanishes on the complement of $K$.    Suppose that $E_J-K$ contains
a point $z_0$.   Then $u=\max_{j\in J}u_j$ in a neighborhood of $z_0$.
Replacing $u$ by
$u^{\epsilon}:=\max(\epsilon_0,u_1+\epsilon_1,\dots,u_A+\epsilon_A)$
for small $\epsilon_j$'s we may assume that $E_J^\epsilon=E_J(\epsilon_0,u_1+\epsilon_1,\ldots,u_A+\epsilon_A)$ is smooth
near $z_0$.  We evaluate $(dd^cu)^{\rm dim}$.
By the useful fact above, near $z_0$ we have $(dd^cu^\epsilon)^N=0$ on $E^{\epsilon}_{\tilde J}- K$ for all ${\tilde J}\subset J$ if $N\ge \sum_{j\in \tilde J}N_j$; thus $(dd^cu^\epsilon)^{\rm dim}=0$.
Letting $\epsilon\to 0$, we have $(dd^cu)^{\rm dim}=0$ near $z_0$, and hence on ${\bf
C}^{\rm dim}-K$.  The formula for $\mu_K$ then follows from Theorem 1.

\bigskip

Let  ${\cal D}^k$ denote the space of test forms of degree $k$.  The
currents of dimension $k$ are defined as the dual of ${\cal D}^k$.
Since we may decompose the $k$-forms into terms of bidegree $(p,q)$,
${\cal D}^k=\bigoplus_{p+q=k}{\cal D}^{p,q}$,  each current may be
written as a sum of currents of bidimension $(p,q)$.  We have
operators $\partial:{\cal D}^{p,q}\to{\cal D}^{p+1,q}$ and
$\bar\partial:{\cal D}^{p,q}\to{\cal D}^{p,q+1}$; and  their
adjoints, which we denote again by $\partial$ and $\bar\partial$,
act on the spaces of currents by duality.  If $u$ is psh, then
$u[E]$ is a current which has the same dimension as $E$, and
$dd^c(u[E])$ is a positive, closed current.  If $u$ is psh and
continuous, we may define $(dd^cu)^n$ by induction on $n$ (cf.
[BT]).  Specifically, since $(dd^cu)^n$ is positive, then it is
represented by integration.    It follows that $u(dd^cu)^n$ is a
well-defined current, and we set
$(dd^cu)^{n+1}:=dd^c(u\,(dd^cu)^n)$, or in other words, its action
on a test form $\varphi$ is given by
$$\langle (dd^cu)^{n+1},\varphi\rangle:= \int dd^c\varphi\wedge u\wedge (dd^cu)^n.$$
This definition gives a continuous extension of $(dd^c)^n$ to the continuous, psh functions.

Let $M\subset{\bf C}^{\rm dim}$ be a smooth submanifold of locally
finite volume.  If $M$ has codimension $k$, then we may orient $M$
by choosing a simple $k$-form $\nu$ of unit length which annihilates
the tangent space to $M$.  We may define the current of integration
$[M]$, which acts on a test form $\varphi$ according to the formula
$$\langle\varphi,[M]\rangle:= \int *(\varphi\wedge\nu)\,||\nu||^{-1}\,dS_M,$$
where $*$ is the Hodge $*$-operator taking volume form to a scalar
function, and $dS_M$ the euclidean surface measure on $M$.  Given a
$k$-tuple $(\rho_1,\dots,\rho_k)$ of defining functions, we define
an orientation as follows.   ${\bf C}^{\rm dim}$ has a canonical
orientation induced by its complex structure.  We orient
$M_1=\{\rho_1=0\}$ as the boundary of $\{\rho_1<0\}\subset{\bf
C}^{\rm dim}$.  Thus $\nu_1=d\rho_1$.  We orient
$M_2=\{\rho_1=\rho_2=0\}$ as the boundary of $\{\rho_2<0\}\cap M_2$
inside $M_1$.  Thus $\nu_2=d\rho_1\wedge d\rho_2$.  Continuing this
way, we orient $M$ using  $\nu=d\rho_1\wedge\cdots\wedge d\rho_k$.

 A $(p,p)$ current $T$ is said to be positive if $\langle T,i\alpha_1\wedge\bar\alpha_1\wedge\cdots\wedge i\alpha_p\wedge\bar\alpha_p\rangle\ge0$ for all smooth test forms $\alpha_j$ of type (1,0).  Here we choose to orient $E_J$ so that $\delta^c_J\wedge[E_J]$ is positive, a choice which is justified by the following.
\proclaim Lemma 1.  We may orient $E_J$ so that $\delta^c_J\wedge
[E_J]$ is a positive current.  If we let $E_J'$ denote $E_J$ with
the orientation given (as above) by taking successive boundaries in
terms of the defining functions
$(\rho_1=u_{j_1}-u_{j_2},\dots,\rho_{\ell-1}=u_{j_{\ell-1}}-u_{j_\ell})$, then we
have $[E_J]=(-1)^{ {\ell(\ell-1)/2}}[E'_J]$.

\noindent{\it Proof. }  The $k$-form defining the orientation of
$E_J$ will be given by $\pm d\rho_1\wedge\cdots\wedge
d\rho_{\ell-1}$, with the sign $\pm$ to be determined.  In the
notation above, we have
$$\delta^c_J\wedge [E_J']=d^c\rho_1\wedge\cdots\wedge d^c\rho_{\ell-1}\wedge d\rho_1\wedge\cdots\wedge d\rho_{\ell-1}\,||d\rho_1\wedge\cdots\wedge d\rho_{\ell-1}||^{-1}\, dS,$$
where $dS$ is the euclidean surface measure on $E_J$. Now since
$d\rho=\partial\rho+\bar\partial\rho$ and
$d^c\rho=i(\bar\partial\rho-\partial\rho)$ we see that
$$(-1)^{  {(\ell-2)(\ell-1)\over 2 }}2^{1-\ell} d^c\rho_1\wedge\cdots\wedge d^c\rho_{\ell-1}\wedge d\rho_1\wedge\cdots\wedge d\rho_{\ell-1}=(i\partial\rho_1\wedge\bar\partial\rho_1)\wedge\cdots\wedge(i\partial\rho_{\ell-1}\wedge\bar\partial\rho_{\ell-1}).$$
It follows that we may orient $E_J$ to make $\delta^c_J\wedge[E_J]$
positive, and the relation between the orientations on $E_J$ and $E'_J$
is as claimed.

\bigskip

We will use the notation:
$$d^c_J=d^c_J(u_1,\dots,u_m):=d^cu_{j_1}\wedge\cdots\wedge d^cu_{j_\ell}. \eqno(4)$$
For $1\leq t\leq\ell$, we write $J(\hat t)$ to denote the set $J$
with the $t$-th element removed; or, if $s\in J$, $J(\hat s)$
denotes the set $J$ with $s$ removed.  The meaning will be clear
from the context.

 \proclaim Lemma 2.   We have the
following identities:
\item{(1)} $ d(d^c_J(u))=\sum_{t=1}^\ell (-1)^{t-1}d^c_{J(\hat t)}(u)\wedge dd^cu_{j_t}.$
\item{(2)} $d^cu_{j_t}\wedge\delta^c_J(u)=(-1)^{\ell-1}d^c_J(u)$ for any $1\le t\le\ell$.
\item{(3)} $ \sum_{t=1}^\ell (-1)^{\ell-t}d^c_{J(\hat t)} = \delta^c_J$.
\item{(4)} $dd^cu_{j_t}\wedge\sigma^n_J+\sigma^{n+1}_{J(\hat t)}=\sigma^{n+1}_J$.

\noindent{\it Proof. }  These identities follow from the product rule and anti-commutation of 1-forms.

\proclaim Lemma 3.  If $\alpha$ is a smooth form of type $(a,a)$,
and if $\beta$ is a smooth form of type $(b,b)$, then the forms
$d\alpha\wedge d^c\beta$ and $d\beta\wedge d^c\alpha$ have the same
parts of type $(a+b+1,a+b+1)$.

\noindent{\it Proof. }  Expand $d\alpha\wedge d^c\beta$ into terms
of the form $\partial\alpha\wedge\bar\partial\beta$, etc., and
compare bidegrees.

\proclaim Lemma 4.  $d[E_J] = [\partial E_J] = \sum_{\tilde J}
\epsilon^J_{\tilde J}[E_{\tilde J}]$, where the sum is taken over
all $\tilde J$ such that $J\subset\tilde J\subset\{1,\dots,m\}$ and
$|\tilde J|=|J|+1$. For each such $\tilde J$, there is an $s$ such
that $j_1<\dots<j_k<s<j_{k+1}<\dots<j_\ell$ and $J=\tilde J(\hat
s)$, and we have $\epsilon^J_{\tilde J}=(-1)^{k}$.

\noindent{\it Proof.}  By Stokes' Theorem, we have $d[E_J]=[\partial
E_J]=\sum [E'_{\tilde J}]$, where $E'_{\tilde J}$ denotes the
manifold $E_{\tilde J}$ with the induced boundary orientation on
$\partial E_J$. Thus we need to compare the orientations of
$E'_{\tilde J}$ and $E_{\tilde J}$.  As in the discussion before
Lemma 1, the orientations of $E_{J}$ and $E'_{\tilde J}$ are given
by the defining functions
$(u_{j_2}-u_{j_1},\dots,u_{j_\ell}-u_{j_{\ell-1}})$ and
$(u_{j_2}-u_{j_1},\dots,u_{j_\ell}-u_{j_{\ell-1}},\rho)$,
respectively, where we may take $\rho$ to be either
$u_s-u_{j_{k+1}}$ or $u_s-u_{j_k}$.  By Lemma 1, the orientation of
$E_J$ is  given by
$$\nu_{\tilde J}=(-1)^{\ell(\ell-1)/2}A\wedge d(u_{j_{k+1}}-u_{j_k})\wedge B,$$
where $A=d(u_{j_2}-u_{j_1})\wedge\cdots\wedge
d(u_{j_k}-u_{j_{k-1}})$ and $B=d(u_{j_{k+2}} - u_{j_{k+1}})\wedge
\cdots\wedge d(u_{j_\ell}-u_{j_{\ell-1}})$.  Thus the orientation of
$E'_{\tilde J}$ is given by the form
$$\nu_{\tilde J}'=(-1)^{\ell(\ell-1)/2}A\wedge d(u_{j_{k+1}}-u_{j_k})\wedge B\wedge d(u_s-u_{j_k}).$$
  By Lemma 1 again, the orientation of $E_{\tilde J}$ is given by the form
$$\nu_{\tilde J}= (-1)^{(\ell+1)\ell/2}A\wedge d(u_s - u_{j_{k}})\wedge d(u_{j_{k+1}}-u_s) \wedge B.$$
Since the degree of $B$ is $\ell-(k+1)$, we find that $\nu'_{\tilde
J} = (-1)^{k}\nu_{\tilde J}$, which completes the proof.

\proclaim Proposition 1.  Let $M$ be a smooth submanifold with
boundary, and let $\chi$ be a smooth form on $E$ so that
$\chi\wedge[M]$ is a current of bidimension $(p,p)$.  If $v$ is a
smooth function, and $\phi$ is a smooth form of bidegree
$(p-1,p-1)$, then
$$\int_M v\wedge \chi\wedge dd^c\phi=\int_M d(d^cv\wedge \chi)\wedge\phi - \int_{\partial M}d^c v\wedge \chi\wedge\phi+\int_M\chi\wedge d(v\wedge d^c\phi).$$

\noindent{\it Proof. }  By the product rule, we have
$$\int_M v\wedge\chi\wedge dd^c\phi=\int_M\chi\wedge\left(d(v\wedge d^c\phi)-dv\wedge d^c\phi\right).$$
Note that $v$ is a (0,0)-form, and $\phi$ is of bidegree
$(p-1,p-1)$, and $\chi\wedge[M]$ is a current of bidimension
$(p,p)$.  Thus the only nonzero terms integrated against  this
current can come from $(p,p)$-forms.  Thus by Lemma 3 we may replace
$dv\wedge d^c\phi$ by $d\phi\wedge d^cv$ in the right hand integral.
Now we integrate by parts in the right hand integral to obtain the
desired formula.
\bigskip

\proclaim Proposition 2.  We have
$$dd^c(u\wedge\delta^c_J\wedge[E_J]) =(-1)^{\ell+1
}\left( d(d^c_J(u))\wedge [E_J]-d^c_J(u)\wedge[\partial E_J] \right)
-d^c(u\wedge d(\delta^c_J\wedge[E_J])).$$

\noindent{\it Proof. }  Let us note first that
$d(\delta^c_J\wedge[E_J])=d(\delta^c_J)\wedge[E_J] + (-1)^{\ell-1}
\delta^c_J\wedge[\partial E_J]$, so this current is represented by
integration.  Thus, since $u=u_j$ for all $j$ on the sets $E_J$ and
$\partial E_J$, the current $u_j\wedge
d(\delta^c_J\wedge[E_J])=u\wedge d(\delta^c_J\wedge[E_J])$ is the
same for all $j\in J$.  Similarly, we may substitute $u_j$ for $u$
in the left hand term of the equation.

By Lemma 1, $\delta^c_J\wedge[E_J]$ is a current of bidegree
$(|J|,|J|)$.  Thus we evaluate the left hand term by testing it
against a form of type $({\rm dim}-|J|-1,{\rm
dim}-|J|-1)$$$\langle\phi,
dd^c(u\wedge\delta^c_J\wedge[E_J])\rangle=\langle
dd^c\phi,u\wedge\delta^c_J\wedge [E_j]\rangle.$$ Now we apply
Proposition 1 with $v=u_j$, $\chi=\delta^c_J$ and  $M=E_J$ to obtain
$$\int_{E_J} u_j\wedge\delta^c_J\wedge dd^c\phi=\int_{E_J}d(d^cu_j\wedge\delta^c_J)\wedge\phi-\int_{\partial E_J} d^cu_j\wedge\delta^c_J\wedge\phi + \int_{E_J}\delta^c_J\wedge d(u_j\wedge d^c\phi).$$
Our formula now follows by applying (2) of Lemma 2 to the first and
second integrals, rewriting the terms as currents, then substituting
$u$ for $u_j$.
\bigskip
\noindent{\bf Proof of Theorem 1. }  In the case $n=1$ we have
$J=\{j\}$, and $\delta^c_J=1$, so the statement of Proposition 2
becomes
$$dd^c(u\wedge [E_j])=dd^c(u_j\wedge[E_j])=dd^cu_j\wedge[E_j]-d^cu_j\wedge[\partial E_j] - d^c(u\wedge d[E_j]).$$
We identify $dd^cu$ with $dd^cu\wedge[E]=dd^c(\sum[E_j])$, which
gives
$$dd^cu\wedge[E]=\sum_j dd^cu_j\wedge [E_j] +\sum_{j_1<j_2}\delta^c_{j_1,j_2}\wedge[E_{j_1,j_2}],$$
where the second sum is a consequence of Lemmas 4 and 2(3) applied
to $\sum d^cu_j\wedge[\partial E_j]$, and the other terms vanish
since $d(\sum[E_j])=d[E]=0$.

Now we proceed by induction, assuming that Theorem 1 has been proved for $n$.  Then we have
$$(dd^cu)^{n+1}=    dd^c(u\wedge(dd^cu)^n)
=   dd^c\left(u\sum_J\sigma_{J}^{n-|J|+1}\wedge\delta^c_J\wedge[E_J]\right). $$
Since $\sigma_{J}^{n-|J|+1}$ is an even form which is both $d$- and $d^c$-closed, this expression is
$$= \sum_J \sigma_{J}^{n-|J|+1}\wedge dd^c(u\wedge\delta^c_J\wedge [E_J]).$$
We apply Proposition 2 to obtain
$$\eqalign{=\sum_J\sigma_{J}^{n-|J|+1} \wedge&\left((-1)^{|J|+1}\big\{d(d^c_J(u)\wedge[E_J] - d^c_J\wedge[\partial E_J]\big\} \right.\cr
&\left. \phantom{(-1)^{|J|+1} }- \big\{d^c(u\wedge d(\delta^c_J\wedge[E_J])\big\}\right)= I + II.\cr}$$
Since $d^c\sigma_J^{n-|J|+1}=0$, we have
$$II=-d^c\sum ud\left(\sigma^{n-|J|+1}\wedge\delta^c_J\wedge[E_J]\right).$$
By our induction hypothesis, then,
$$II= - d^c(u\wedge d((dd^cu)^n)= -d^c(u\wedge 0)=0.$$

Now use (1) of Lemma 2 in the left hand summation in $I$ to obtain
$$\eqalign{(dd^cu)^{n+1}&=\sum\sigma^{n-|J|+1}_J\wedge (-1)^{|J|+1}\sum_{t=1}^{|J|}(-1)^{t-1}dd^cu_{j_t}\wedge d^c_{J(\hat t)}\wedge[E_J]\cr  &\ -\ \sum(-1)^{|J|+1} \sigma_J^{n-|J|+1}\wedge d^c_J\wedge[\partial E_J]
=A+B.}$$
In the notation of Lemma 4, we have
$$B=\sum_J \sigma^{n-|J|+1}_J\wedge(-1)^{|J|} d^c_J\wedge\sum_{\tilde J}\epsilon^J_{\tilde J}[E_{\tilde J}].$$
Now let us rewrite $B$, summing over $\tilde J$ on the outside, and
summing over subsets $J=\tilde J(\hat s)$ on the inside.  By Lemma 4,
 $\epsilon_{\tilde J}^{\tilde J(\hat s)}= (-1)^{s-1}$, $1\leq s\leq |\tilde J|$.  This gives
$$B= \sum_{\tilde J}\left (\sum_{s=1}^{|\tilde J|} (-1)^{|\tilde J|-s} \sigma^{n-|\tilde J|+2}_{\tilde J(\hat s)}\wedge d^c_{\tilde J(\hat s)} \right)\wedge[E_{\tilde J}].$$
Since we are summing over all subsets $\tilde J$, we can remove the
tilde from $\tilde J$.  Further, we can set $s=t$, which lets us
rewrite $A+B$ as
$$\sum_{J} \sum_t\left( \sigma^{n-|J|+1}_J\wedge dd^cu_{j_t} + \sigma^{n-|J|+2}_{J(\hat t)}\right)\wedge d^c_{J(\hat t)}\wedge(-1)^{|J|-t} d^c_J(\hat t)\wedge[E_J].$$
Finally, by (3) and (4) of Lemma 2, we have
$$(dd^cu)^{n+1}=\sum_J \sigma^{n-|J|+2}_J\wedge\left(\sum_t (-1)^{|J|-t} d^c_{J(\hat t)}\right)\wedge [E_J] =
\sum_J \sigma^{n-|J|+2}_J\wedge \delta^c_J\wedge [E_J]$$
which completes the proof.

\bigskip

\noindent {\bf Acknowledgement.} The authors would like to thank
Norm Levenberg for helpful discussions on this material.

\bigskip
\bigskip
\centerline{\bf References}
\bigskip

\item{[BT]}  E. Bedford and B.A. Taylor, The Dirichlet problem for a complex Monge-Amp\`ere equation, {\it Invent. Math.} {\bf 37} (1976), no. 1, 1--44.

\item{[B1]}  Z. B\l ocki,  Equilibrium measure of a product set of ${\bf C}^n$, {\it Proc. A.M.S.} {\bf 128} (2000), no.\ 12, 3595--3599.

\item{[B2]}  Z. B\l ocki, The domain of definition of the complex Monge-Amp\`ere operator, {\it Amer. J. Math.} {\bf 128} (2006), no. 2, 519--530.

\item{[M]}  S. Ma`u, {\it Plurisubharmonic Functions of Logarithmic
Growth}, PhD Thesis, University of Auckland, 2003.

\bigskip

\rightline{ Indiana University}

\rightline{Bloomington, IN 47405 USA}
\smallskip
\rightline{\tt bedford@indiana.edu}

\rightline {\tt sinmau@indiana.edu}

 \bye